\documentclass[11pt]{amsart}
\usepackage{amssymb}
\theoremstyle{plain}
\newtheorem{theorem}{Theorem}[section]
\newtheorem{corollary}[theorem]{Corollary}
\newtheorem{lemma}[theorem]{Lemma}

\newtheorem{fact}[theorem]{Fact}
\newtheorem{claim}[theorem]{Claim}

\theoremstyle{definition}
\newtheorem{definition}[theorem]{Definition}
\newtheorem{remark}[theorem]{Remark}

\theoremstyle{remark}

\newcommand{\id}{\operatorname{id}}
\newcommand{\LS}{\operatorname{LS}}

\newcommand{\cf}{\operatorname{cf}}

\newcommand{\gaS}{\operatorname{ga-S}}
\newcommand{\Hanf}{\operatorname{Hanf}}

\renewcommand{\phi}{\varphi}

\newcommand{\Union}{\bigcup}
\newcommand{\initial}\lessdot
 
\newcommand{\K}{\operatorname{\mathcal{K}}}

\newcommand{\C}{\mathfrak C}

\def\?{?\vadjust
{\vbox to 0pt{\vskip-7pt\hbox to 1.1\hsize{\hfill\huge ?!}}}}

\def\bK{\K}
\def\subm{\prec_{\bK}}
\begin{document}
\title{Upward Stability Transfer for Tame Abstract Elementary Classes}
\author{John Baldwin}
\email[John Baldwin]{jbaldwin@uic.edu}
\address{Department of Mathematics\\
University of Illinois at Chicago\\
Chicago IL 60607}

\author{David Kueker}
\email[David Kueker]{dwk@math.umd.edu}
\address{Department of Mathematics\\
University of Maryland\\
College Park MD 20742-4015}

\author{Monica VanDieren}
\email[Monica VanDieren]{mvd@umich.edu}
\address{Department of Mathematics\\
University of Michigan\\
Ann Arbor MI 48109-1043}


%\date{February 23, 2004}
%\date{September 15, 2005}
\date{October 29, 2005}

\begin{abstract}
Grossberg and VanDieren have started a program to develop a
stability theory for tame classes (see \cite{GrVa}).  We name
some variants of tameness (Definitions \ref{tamedef} and
\ref{wtamedef}) and prove the following.

\begin{theorem}\label{up stab}
Let $\K$ be an AEC with L\"{o}wenheim-Skolem number $\leq\kappa$.
Assume that $\K$ satisfies the amalgamation property and is
$\kappa$-weakly tame and Galois-stable in $\kappa$.  Then, $\K$ is
Galois-stable in $\kappa^{+n}$ for all $n<\omega$.
\end{theorem}

With one further hypothesis we get a very strong conclusion in the
countable case.

\begin{theorem}\label{genbasecase} Let $\K$ be an AEC satisfying
the amalgamation property and with L\"{o}wenheim-Skolem number
$\aleph_0$ that is $\omega$-local and $\aleph_0$-tame. If $\K$ is
$\aleph_0$-Galois-stable then $\bK$ is Galois-stable in all
cardinalities.
\end{theorem}
\end{abstract}

\maketitle

%%%%%%%%%%%%%%%%%%%%%%%%%%%%%%%
\bigskip
\section*{Introduction} \label{s:introductionm}

A tame abstract elementary class  is an abstract elementary class
(AEC) in which inequality of Galois-types has a local behavior.
Tameness is a natural condition, generalizing both homogeneous
classes and excellent classes, that has very strong consequences.
We examine one of them here.

The work discussed in this paper fits in the program of developing a
model theory, in particular a stability theory, for non-elementary
classes.  Many results to this end were in contexts
where manipulations with first order formulas, or infinitary formulas,
were pertinent and consequential.  Most often, types in these context were
identified with satisfiable collections of formulas.  The model theory for
abstract elementary classes where types are identified roughly with the
orbits of an element under automorphisms of some large structure moves
away from the dependence on ideas from first order logic.

The main result of this paper is not surprising in light of what
is known about first order model theory, but it does shed light on
problems that become more  elusive  in abstract elementary
classes. Grossberg and VanDieren \cite{GrVa} provide a sufficient
condition for stability which yields Theorem 0.1 under GCH. They
prove (Their paper assumes $\mu$ greater than the Hanf number but
this is not needed; see \cite{Baldwin_monograph}.):

\begin{fact}[Corollary 6.4 of \cite{GrVa}]

Suppose that $\K$ is a $\chi$-tame AEC for some $\chi\geq
\LS(\bK)$. If $\K$ is Galois-stable in some $\mu\geq \chi$,
 then $\K$ is stable in every
$\kappa$ with $\kappa^\mu =\kappa$.

\end{fact}

The \cite{GrVa} argument
generalizes an aspect of Shelah's method for calculating the
entire spectrum function. The ZFC-argument here illustrates the
relation of tame AEC's to first order logic.  We adapt an argument
for Morley's Theorem that $ \omega$-stability implies stability in
all cardinalities to the context of Galois-types to move
Galois-stabililty from a cardinal to its successor. For larger
$\kappa$  some splitting technology is needed and the result is
that Galois-stability in $\kappa$ implies Galois-stability in
$\kappa^+$ when $\kappa$ is at least as large as the tameness
cardinal.

Combining the result of \cite{GrVa} along with the results of this paper,
we gain a better understanding of the stability spectrum for tame AECs:

\begin{corollary}[Partial Stability Spectrum]

Suppose that $\K$ is a $\chi$-tame AEC for some $\chi\geq
\LS(\bK)$. If $\K$ is Galois-stable in some $\mu \geq \chi$, then
$\K$ is stable in every $\kappa$ with $\kappa^\mu =\kappa$ and in
$\mu^{+n}$ for all $n<\omega$.

\end{corollary}

We acknowledge helpful conversations with Rami Grossberg and Alexei
Kolesnikov, particularly on the correct formulation and proof of Fact
\ref{exist univ ext}.

\section{Background}
Much of the necessary background for this paper can be found in
the exposition \cite{Gr1} and the following papers on tame abstract
elementary classes \cite{GrVa} and \cite{GrVa2}.
We will review some of the required definitions and theorems in this
section.
We will use
$\alpha,\beta,\gamma,i,j$ to denote ordinals and
$\kappa,\lambda,\mu,\chi$ will be used for cardinals.
We will use $(\K,\prec_{\K})$ to denote an abstract
elementary class  and $\K_\mu$ is the subclass of models in $\K$ of
cardinality $\mu$.
For an AEC $\K$, $\LS(\K)$ represents the L\"{o}wenheim-Skolem number of
the class. Models are denoted by $M,N$ and may be
decorated with superscripts and subscripts.  Sequences of elements
from $M$ are written as $\bar a,\bar b,\bar c,\bar d$. The letters
$e,f,g,h$ are reserved for $\K$-mappings and $\id$ is the identity
mapping.

For the remainder of this paper we will fix $(\K,\prec_{\K})$ to
be an abstract elementary class {\em satisfying the amalgamation
property}.  It is easy to see that we only make use of the
$\kappa$-amalgamation property for certain $\kappa$ and some facts
here hold in classes satisfying even weaker amalgamation
hypotheses. Since we assume the amalgamation property, we can fix
a monster model $\C\in\mathcal K$ and  say that the type of $a$
over a model $M\prec_{\K}\C$ is equal to the type $b$ over $M$ iff
there is an automorphism of $\C$ fixing $M$ which takes $a$ to
$b$. (Technically, the existence of a monster model requires the
joint embedding property as well as the amalgamation property.
However, in the presence of the amalgamation property, joint
embedding is an equivalence relation and our fixing of the monster
model is the same as restricting to one equivalence class.)  In
this paper, we work entirely with \emph{Galois-types} (i.e.
orbits) and so feel free to write simply type.
% In this paper we
%will freely use the term \emph{type} in place of
%\emph{Galois-type} which is used in the literature to distinguish
%types defined by collections of formulas from those defined as
%orbits.
For a model $M$ in $\K$, the set of Galois-types over $M$
is written as $\gaS(M)$. An AEC $\K$ satisfying the amalgamation
property is \emph{Galois-stable} in $\kappa$ provided that for
every $M\in\K_\kappa$ the number of types over $M$ is
$\leq\kappa$.

Let us recall a few results that follow from
Galois-stability in
$\kappa$.

\begin{definition}
Let $M\in\K_\kappa$, we say that \emph{$N$ is universal over $M$}
provided that
for every $M'\in\K_{\kappa}$ with $M\prec_{\K}M'$, there exists a
$\K$-mapping $f:M'\rightarrow N$ such that $f\restriction M=id_M$.
\end{definition}

Note that in contrast to most model theoretic literature, in AECs
a tradition has grown up of defining `universal over $M$' as
`universal over submodels of the same size as $M$'.

\begin{fact}[\cite{Sh 600}, see  \cite{Baldwin_monograph} or \cite{GrVa}
for a proof]\label{exist univ ext}\label{getu}
If $\K$ is Galois-stable in $\kappa$ and satisfies the
$\leq\kappa$-amalgamation property, then for every
$M\in\K_{\kappa}$ there is some (not necessarily unique) extension $N$ of
$M$ of cardinality
$\kappa$ such that $N$ is universal over $M$.
\end{fact}

If $\K$ is Galois-stable in $\kappa$, we can construct an
increasing and continuous chain of models $\langle
M_i\in\K_{\kappa}\mid i<\sigma\rangle$ for any limit ordinal
$\sigma\leq\kappa^+$ such that $M_{i+1}$ is universal over $M_i$.
The limit of such a chain is referred to as a
{\emph{$(\kappa,\sigma)$-limit model}.

\begin{corollary}\label{getsat}  Suppose $\bK$ is $\kappa$-Galois-stable and
$\bK_{\kappa}$ has
the amalgamation property with $\LS(\K)\leq\kappa$.
Then for any model $M \in \bK$ with
cardinality $\kappa^+$ we can find a $\kappa^+$-saturated and
$(\kappa,\kappa^+)$-limit model $M'$ such that $M$ can be embedded
in $M'$.
\end{corollary}

\begin{proof} Write $M$ as an increasing continuous chain $M_i$ of
models of cardinality at most $\kappa$.  We define an increasing
chain of  models $M'_i$, each with cardinality $\kappa$, and $f_i$
so that $f_i$ is a $\bK$-embedding of $M_i$ in $M'_i$ and such
that each $M'_{i+1}$ realizes all types over $M_{i}$; indeed,
$M'_{i+1}$ is universal over $M'_{i}$. For this, first choose
$M^1_i$ which is universal over $M'_{i}$ by Fact \ref{exist univ
ext}. Then amalgamate $M_{i+1}$ and $M^1_i$ over $f_i:M_i \mapsto
M'_i$ with $M'_i\subm M'_{i+1}$. Now the union of the $M'_i$ is a
$(\kappa,\kappa^+)$-limit model which imbeds $M$.
\end{proof}

Now we turn our attention to two definitions which capture instances in
which types are determined by a small set.  These two approaches to local
character play different roles in this paper.

\begin{definition}\label{tamedef}
 Let $\K$ be an AEC.
\begin{enumerate}
\item We  say
that a class $\K$ is \emph{$\chi$-tame} provided that for every
model $M$ in $\K$  with $|M| \geq \chi$ and every $p$ and $q$,
types over $M$, if $p\neq q$, then there is a model of cardinality
$\chi$ which distinguishes them.  In other words if $p\neq q$,
then there exists $N\in\K_\chi$ with $N\prec_{\K}M$ such that
$p\restriction N\neq q\restriction N$.

\item A class $\K$ is \emph{$\omega$-local} provided for every increasing
chain of
types $\{p_i\mid i<\omega\}$ there is a unique $p$ such that
$p=\Union_{i<\omega}p_i$.

\end{enumerate}
\end{definition}

 For some of the results in this paper we could
replace $\chi$-tameness with the two-parameter version of
\cite{BaldwinEM}, $(\kappa,\chi)$-tameness, which requires  only
that distinct types over models of cardinality $\kappa$ be
distinguished by models of cardinality $\chi$.  Since we don't
actually carry out any inductions to establish tameness, this
nicety is not needed here.  Note that if $\chi< \kappa$,
$\chi$-tame implies $\kappa$-tame.

\begin{remark}\label{union of types}
If $\K$ is an AEC with the amalgamation property, for every
increasing $\omega$-chain of types $p_i$, there is a type over the
union of the domains extending each of the $p_i$ (1.10 of \cite{Sh
394}, proved as 3.14 in \cite{BaldwinEM}); however, this extension
need not be unique.
\end{remark}

\begin{remark}
%If  an AEC axiomatized by a $L_{\omega_1,\omega}(\mathbf
%Q)$-sentence satisfies \\
Clearly, if an AEC is defined by a logic with finitary syntax, has
L\"owenheim-Skolem number $\aleph_0$, and  ``Galois-types =
syntactic types" then the AEC is both $\aleph_0$-tame and
$\omega$-local. Shelah showed, assuming weak GCH,  this  happens
for $L_{\omega_1,\omega}$ classes that are categorical in
$\aleph_{n}$ for every $n<\omega$; it also holds for Zilber's
quasiminimal excellent classes.
% $L_{\omega_1,\omega}(\mathbf Q)$.
\end{remark}

 A
weaker version of tameness requires that only those types over
saturated models are determined by small sets. This appears as
$\chi$-character in \cite{Sh 394} where Shelah proves that, in
certain situations, categorical AECs have small character.

\begin{definition}\label{wtamedef}
For an AEC $\K$ and a cardinal $\chi$, we say that $\K$ is
\emph{$\chi$-weakly tame} or has \emph{$\chi$-character} iff for
every saturated model $M$ with $|M| \geq \chi$ and every $p\neq
q\in \gaS(M)$, there exists $N\in\K_\chi$ such that $N\prec_{\K}M$
and $p\restriction N\neq q\restriction N$.
\end{definition}

%%%%%%%%%%%%%%%%%%%%%%%%%%%%%
\bigskip
\section{$\aleph_0$-tameness} \label{s:omega stability}
In this section we assume $\bK$ has a countable language,
L\"{o}wenheim-Skolem number $\omega$and is $\aleph_0$-tame.
% since
%in this case we get a more general result than stated in the
%abstract.

\begin{theorem}\label{basecase} Suppose $\LS(\K)=\aleph_0$.
If $\K$ is $\aleph_0$-tame and $\mu$-Galois-stable for all $\mu <
\kappa$ and $\cf(\kappa) > \aleph_0$ then $\bK$ is
$\kappa$-Galois-stable.
\end{theorem}

\begin{proof}
For purposes of contradiction suppose there are more than $\kappa$
types over some model $M^*$ in $\K$ of cardinality $\kappa$.
 We
may  write $M^*$ as the union of a continuous chain $\langle
M_i\mid i<\kappa\rangle$ under $\prec_{\K}$ of models in $\K$
which have  cardinality $<\kappa$. We say that a type over $M_i$
\emph{has many extensions} to mean that it has $\geq \kappa^+$
distinct extensions to a type over $M^*$.

\begin{claim}\label{some type has many ext}
 For every $i$, there is some type over $M_i$ with many
extensions.
\end{claim}

\begin{proof}[Proof of Claim \ref{some type has many ext}]
    Each type over $M^*$ is the extension of some type over $M_i$ and, by
our assumption, there are less than $\kappa$ many types over
$M_i$, so at least one of them must have many extensions.
%\renewcommand{\qedsymbol}{\dashv_{\ref{some type has many ext}}}
%${_{\ref{some type has many ext}}}$

\end{proof}%${_{\ref{some type has many ext}}}$

\begin{claim}\label{p with many has ext with many}
 For every $i$, if the type $p$ over $M_i$ has many extensions,
then for every $j>i$, $p $ has an extension to a type  $p'$ over $M_j$
with many extensions.
\end{claim}

\begin{proof}[Proof of Claim \ref{p with many has ext with many}]
  Every extension of $p$ to a type over $M^*$ is the extension of some
extension of $p$ to a type over $M_j$. By our assumption there are
less than $\kappa$ many such extensions to a type over $M_j$, so
at least one of them must have many extensions.
\end{proof}%$_{\ref{p with many has ext with many}}$

\begin{claim}\label{two ext}
For every $i$, if
 the type $p$ over $M_i$ has many
extensions, then for all sufficiently large $j>i$, $p$ can be
extended to two types over $M_j$ each having many extensions.
\end{claim}

\begin{proof}[Proof of Claim \ref{two ext}]
By Claim \ref{p with many has ext with many} it suffices to
establish the result for some $j>i$. So assume that there is no
$j>i$ such that $p$ has two extensions to types over $M_j$ each
having many extensions. Then, by Claim \ref{p with many has ext
with many} again, for every $j>i$, $p$ has a unique extension to a
type $p_j$ over $M_j$ with many extensions. Let $S^*$ be the set
of all extensions of $p$ to a type over $M^*$ -- so
$|S^*|\geq\kappa^{+}$. Then $S^*$ is the union of $S_0$ and $S_1$,
where  $S_0$ is the set of all $q$ in $S^*$ such that $p_j
\subseteq q$ for all $j>i$,
 and $S_1$ is the set of all $q$ in $S^*$ such
that $q$ does not extend $p_j$ for some $j>i$. Now if $q_1$ and
$q_2$ are different types in $S^*$ then, since $\K$ is
$\aleph_0$-tame and $\cf(\kappa) > \aleph_0$,
%\sidebar{ OLD:their restrictions to some $M\prec_{\K}M^*$ of cardinality
%$\kappa$ must differ.}
their restrictions to some
$M_i\prec_{\K}M^*$ with $i< \kappa$ must differ. Hence their
restrictions to all sufficiently large $M_j$ must differ.
Therefore, $S_0$ contains at most one type. On the other hand, if
$q$ is in $S_1$ then, for some $j>i$, $q\restriction M_j$ is an
extension of $p$ to a type over $M_j$ which is different from
$p_j$, hence has at most $\kappa$ extensions to a type over $M^*$.
Since there are $<\kappa$ types over each $M_j$ (by our stability
assumption) and just $\kappa$ models $M_j$ there can be at most
$\kappa$ types in $S_1$. Thus $S^*$ contains at most $\kappa$
types, a contradiction.

\end{proof}%${\ref{two ext}}$

\begin{claim}\label{gettree} There is a countable $M \subm
M^*$ such that there are $2^{\aleph_0}$ types over $M$.
\end{claim}

\begin{proof}[Proof of Claim \ref{gettree}]
  Let $p$ be a type over $M_0$ with many extensions.  By
Claim~\ref{two ext} there is a $j_1 > 0$ such that $p$ has two
extensions $p_0, p_1$ to types over $M_{j_1}$ with many
extensions.  Iterating this construction we find a sequence of
 models $M_{j_n}$ and a tree $p_s$ of types for $s \in
2^{\omega}$ with the $2^n$ types $p_s$ (where $s$ has length $n$)
all over $M_{j_n}$ and each $p_s$ has  many extensions. Invoking
$\aleph_0$-tameness, we can replace each $M_{j_n}$ by a countable
$M'_{j_n}$ and $p_s$ by $p'_s$ over $M'_{j_n}$ while preserving
the tree structure on the $p'_s$.  Let $\hat M$ be the union of
the $M'_{j_n}$. Now for each $\sigma \in 2^{\omega}$, $p_\sigma =
\bigcup_{s \subset \sigma} p_s$ is a Galois-type, by Remark
\ref{union of types}
\end{proof}%$\ _{\ref{gettree}}$

Since Claim~\ref{gettree} contradicts the hypothesis of
$\omega$-Galois-stability, this establishes Theorem~\ref{basecase}.
\end{proof}%$\ _{\ref{basecase}}$

Now we obtain Theorem \ref{genbasecase} from the abstract.

\begin{corollary}\label{actgenbasecase}
Suppose $\LS(\K)=\aleph_0$ and $\bK$ has the amalgamation
property. If $\K$ is $\aleph_0$-weakly-tame and
$\omega$-Galois-stable then
\begin{enumerate}
\item\label{actgenbasecase1}
$\bK$ is Galois-stable in all $\aleph_n$ for $n< \omega$.

\item If in addition $\bK$ is both $\omega$-local and $\aleph_0$-tame,
$\bK$ is Galois-stable in
all cardinalities.
\end{enumerate}
\end{corollary}

\begin{proof}[Proof of Corollary \ref{actgenbasecase}]
In the proof of Theorem \ref{basecase} if $\kappa$ is a successor cardinal,
then by Corollary \ref{getsat},
$M^*$ can be embedded into a saturated model and the proof can be carried
through with the weaker assumption of
$\aleph_0$-weak-tameness.
Thus the first claim follows by induction.

To carry out the induction for all cardinals,  we extend the
argument in Theorem \ref{basecase} to limit cardinals of
cofinality $\omega$. At the stage where we called upon
$\aleph_0$-tameness in Claim \ref{two ext}, we now use the
hypothesis of $\omega$-locality. For limit cardinals of
uncountable cofinality, we use the assumption of
$\aleph_0$-tameness since we have no guarantee that $M^*$ can be
taken to be saturated.

\end{proof}

%%%%%%%%%%%%%%%%%%%%%%%%%%%%%%%
\bigskip
\section{$\kappa$-tame: Uncountable $\kappa$} \label{s:kappa stability}

Note that the proof of Theorem \ref{basecase} cannot be immediately
generalized to deducing stability in $\kappa^+$ from stability in
$\kappa$ when the class is tame, but not $\aleph_0$-tame.  The fact that
the countable increasing union of Galois types is a Galois type is very
much particular to `countable' and in general does not hold when we
replace countable by uncountable.    We solve this with a use of
$\mu$-splitting.

\begin{definition}[\cite{Sh
394}]\index{$\mu$-splits}\index{Galois-type!$\mu$-splits} A type
$p\in \gaS(N)$ \emph{$\mu$-splits} over $M\prec_{\K}N$ if and only if
 there exist $N_1,N_2\in\K_{\leq\mu}$ and $h$, a
$\K$-embedding such  that $M\prec_{\K}N_l\prec_{\K}N$ for $l=1,2$ and
$h:N_1\rightarrow N_2$ such that $h\restriction M=\id_M$ and
$p\restriction N_2\neq h(p\restriction N_1)$.

\end{definition}

This dependence relation behaves nicely in Galois-stable AECs. The
existence of unique non-splitting extensions from $M$ to $M'$
where $M$ and $M'$ have the same cardinality and $M'$ is universal
over $M$ holds for any AEC with amalgamation. There is a full
proof as 1.4.13 and 1.4.14 of \cite{vandieren}.  Existence of
non-splitting extensions to larger cardinalities is more
difficult;  under the assumption of categoricity, such an
extension property is asserted in \cite{Sh 394} and a special case
is given a short proof in \cite{Baldwin_nonsplit}.  In the more
general situation, uniqueness requires tameness; see 6.2 of
\cite{Sh 394}. Here we state the uniqueness and existence
statements upon which we will be explicitly calling.

\begin{lemma}[Uniqueness \cite{Sh 394} and \cite{vandieren}]\label{unique
ext}\index{non-$\mu$-splitting!uniqueness} Let $N,M,M'\in\K_\mu$ be
such that
$M'$ is universal over
$M$ and
$M$ is universal over $N$.  If $p\in \gaS(M)$ does not $\mu$-split over
$N$, then there is a unique $p'\in\gaS(M')$ such that $p'$ extends $p$
and $p'$ does not $\mu$ split over $N$.
\end{lemma}

\begin{lemma}[Existence Fact 3.3 of \cite{Sh 394} see also
\cite{GrVa}]\label{exist non-split} Let $M\in\K_{\geq\kappa}$ be
given. Suppose that $\K$ satisfies the $(\leq\|M\|)$-amalgamation
property.  If $\K$ is Galois-stable in $\kappa$, then for every
$p\in \gaS(M)$, there exists $N\in\K_\kappa$ such that
$N\preceq_{\K}M$ and $p$ does not $\kappa$-split over $N$.
\end{lemma}

\begin{remark}  The  arguments in
Claim~\ref{gettree} and Lemma~\ref{exist non-split} differ.  In
Claim~\ref{gettree}, we construct a tree of height $\omega$ of
Galois types and must find a limit for each branch. In
Lemma~\ref{exist non-split}, a tree of height $\kappa$ is
constructed by spreading out copies of a given type.
\end{remark}

 We are able to carry out the following
argument under the hypothesis of weakly tame rather than tame so
we record the stronger result.

\begin{theorem}\label{kappacase}
Let $\K$ be an abstract elementary class with the amalgamation
property that has L\"{o}wenheim-Skolem number $\leq\kappa$ and is
$\kappa$-weakly-tame. Then if $\K$ is Galois-stable in $\kappa$ it
is also Galois-stable in $\kappa^+$.
\end{theorem}

\begin{proof} We proceed by contradiction. So we make the following
assumption: $M^*$ is a model of cardinality $\kappa^+$ with more
than $\kappa^+$ types over it.  By Corollary~\ref{getsat}, we can
extend $M^*$ to a $(\kappa,\kappa^+)$-limit model which is saturated.
Since it has at least as many types as the original we just assume that
$M^*$ is a saturated, $(\kappa,\kappa^+)$-limit model witnessed by $\langle
M_i\mid i<\kappa^+\rangle$.

Let $\{p_\alpha \mid \alpha<\kappa^{++}\}$ be a set of distinct
types over $M^*$.  By stability in $\kappa$, for every $p_\alpha$
there exists $i_\alpha<\kappa^+$ such that $p_\alpha$ does not
$\kappa$-split over $M_{i_\alpha}$ (by
Lemma \ref{exist non-split}). (Note, we don't need a
$(\kappa,\kappa^+)$-limit here but we do below.) By the
pigeon-hole principle there exists $i^*<\kappa^+$ and $A\subseteq
\kappa^{++}$ of cardinality $\kappa^{++}$ such that for every
$\alpha\in A$, $i_\alpha=i^*$.

Now apply the argument of the Claims from the previous section to the
$p_\alpha$ for
$\alpha\in A$ to conclude there exist $p,q \in S(M^*)$ and $i<i'\in A$, such that
neither $p$ nor $q$ $\kappa$-splits over $M_i$ or $M_{i'}$ but
$p\restriction M_{i'} = q \restriction M_{i'}$.
By weak tameness, there exists an ordinal $j>i'$ such that
$p\restriction M_j\neq q\restriction M_j$.  Notice that neither
 $p\restriction
M_j$ nor $q\restriction M_j$ $\kappa$-split over $M_i$.
This
contradicts
Lemma
\ref{unique ext} by giving us two distinct extensions of a non-splitting
type to the model $M_j$ which by construction is universal over $M_{i'}$.
\end{proof}

Using Theorem \ref{kappacase} with an inductive argument on
$n<\omega$, together with the argument for Corollary
\ref{actgenbasecase} (\ref{actgenbasecase1}), we obtain
 Theorem 1 from the abstract:

\begin{theorem}\label{successors}
Let $\K$ be an abstract elementary class that has
L\"{o}wenheim-Skolem number $\leq\kappa$ and satisfies the
amalgamation property and is  $\kappa$-weakly tame. Then if $\K$
is Galois-stable in $\kappa$ it is also Galois-stable in
$\kappa^{+n}$ for any $n<\omega$.
\end{theorem}

One motivation for working out these arguments was to explore
whether or not Galois-superstability (in the sense of few types
over models in every large enough cardinality) could be derived
from categoricity in the abstract elementary class setting.
Following tradition, we write $\Hanf(\K)$ for the Hanf number for
omitting types in first order languages with the same size
vocabulary as $\bK$.  Using Ehrenfeucht-Mostowski models as in the
first order case, for an AEC with amalgamation, categoricity in a
$\lambda$ greater than $\Hanf(\K)$ implies Galois-stability below
$\lambda$.  In the first order case, analysis of the stability
spectrum function allows one to conclude stability in $\lambda$.
Although we don't have such a full analysis of the spectrum
function, we can immediately conclude from
Theorem~\ref{kappacase}:

\begin{corollary}\label{cat and stability}  Suppose
 $\lambda$ is a successor cardinal greater than $\Hanf(\K)$.
  Let $\K$ be an abstract elementary class with the amalgamation
  property
that has L\"{o}wenheim-Skolem number $<\lambda$ and is
$\lambda$-weakly tame. If $\K$ is $\lambda$-categorical, then
 it is  Galois-stable in $\lambda$.
 \end{corollary}

This result is also a
consequence of Theorem 4.1 in \cite{GrVa2} in which the hypotheses of
Corollary \ref{cat and stability} allow one to construct for every
$M\in\K_\lambda$ a model $M'$ also of cardinality $\lambda$ so that $M'$
realizes every type over $M$.

%%%%%%%%%%%%%%%%%%%%%%%%%%%%%%%%%%%%%%%%%%%%%%%%%%%%%%%%%%%%%%%%%%
%%%%%%%%%%%%%%%%%%%%

\end{document}